# A NONMANIPULABLE TEST

By Wojciech Olszewski[1] and Alvaro Sandroni[2]

*Northwestern University and University of Pennsylvania*

A test is said to control for type I error if it is unlikely to reject the data-generating process. However, if it is possible to produce stochastic processes at random such that, for all possible future realizations of the data, the selected process is unlikely to be rejected, then the test is said to be manipulable. So, a manipulable test has essentially no capacity to reject a strategic expert.

Many tests proposed in the existing literature, including calibration tests, control for type I error but are manipulable. We construct a test that controls for type I error and is nonmanipulable.

**1. Introduction.** Professional forecasts are often presented as probabilistic statements [see, e.g., Gneiting and Raftery (2005) for a review of the role of probabilistic forecasts in meteorology]. The quality of these forecasts is regularly tested empirically. A concern [which can be traced back at least to Brier (1950)] is that if forecasts are tested, then experts may misreport their forecasts with the intention of passing the test. Recent literature shows that, without any knowledge about the data-generating process, it is possible to produce forecasts that pass some empirical tests on all possible future realizations of the data.

An example of a test that can be manipulated in this way is the well-known calibration test. Suppose that a stochastic process generates in every period an outcome that can be either 0 or 1. The calibration test requires the empirical frequency of 1 to be close to $p$ in the periods that 1 was forecasted with probability close to $p$. Dawid (1982) shows that the forecasts of the data-generating process will eventually be calibrated. Foster and Vohra (1998) show that any individual can produce forecasts with a random device such that, for all possible infinite strings of zeros and ones, the realized fore-

Received May 2007; revised January 2008.

[1]Supported in part by NSF CAREER Grant SES-0644930.

[2]Supported in part by NSF Grant SES-0820472.

casts will eventually be calibrated as well, with probability one according to the random device used by the individual.

It is natural to seek tests that cannot be manipulated. We address this problem in the present paper. We study the following framework: A stochastic process generates in every period an outcome that can be 0 or 1 (we make no assumptions, such as an independent, identically distributed process). Before any data are observed, an expert named Bob delivers a theory, defined as a mapping that takes as an input any finite string of outcomes and returns as an output a probability of 1; equivalently, theories are probability measures $P$ on the space of infinite sequences of outcomes.

A tester named Alice tests Bob's theory $P$ by selecting an event $A_P$ (i.e., a set of sequences of outcomes), which Alice regards as consistent with the theory $P$. We call $A_P$ the acceptance set and its complement the rejection set for the theory $P$. If, for every $P$, the event $A_P$ has high probability according to $P$, then the data-generating process will not be rejected (with high probability). We then say that the test controls for type I error (of rejecting the data-generating process).

Assume that Bob knows nothing about the data-generating process. However, Bob may use a random device $\zeta$ to select his theory $P$. At first, Alice cannot tell whether the announced theory coincides with the data-generating process or was selected at random. This must be determined by the data. If, for any sequence of outcomes, Bob's theory $P$ is not rejected with high probability according to Bob's random device $\zeta$, then we say that the test can be manipulated (with this high probability), or we say, equivalently, that the test has essentially no capacity to reject theories produced strategically.

Many tests that control for type I error have essentially no capacity to reject theories produced strategically. The calibration test can be manipulated. Several extensions of the calibration test have also been proven to be manipulable [see, e.g., Fudenberg and Levine (1999), Lehrer (2001) and Sandroni, Smorodinsky and Vohra (2003)]. Other statistical tests can also be manipulated. A *prequential test* rejects or accepts a theory based only on the observed data sequence and the next period forecasts made by the theory along the realized sequence of outcomes [see Dawid (1991) for more details on the prequential principle]. Many standard statistical tests, including calibration tests, are prequential tests, and prequential tests can be manipulated [see Sandroni (2003), Vovk and Shafer (2005), Olszewski and Sandroni (2008) and Shmaya (2008)].

Dekel and Feinberg (2006) show that, under the continuum hypothesis, there exists a test that does not reject the data-generating process and cannot be manipulated because every random device $\zeta$ fails this test with certainty on an uncountable number of paths. We construct a (nonprequential) test, called *the global category test*, that also does not reject the data-generating process with probability one and cannot be manipulated. However, our construction does not assume the continuum hypothesis; that is, it is performed within the Zermelo–Fraenkel axioms and the axiom of choice. In addition, the global category test is explicitly constructed; that is, we give set-theoretic formulas for the acceptance sets.

The significance of dispensing the continuum hypothesis can be seen in a simple example. Assume that Bob announces an easy to describe theory (e.g., 0 and 1 occur with probability 0.5 in all periods) and the sequence of outcomes follows an easy-to-describe deterministic process (e.g., 0 in every third period and 1 otherwise). If Alice uses the test in Dekel and Feinberg (2006), no researcher can determine whether Alice rejects Bob's theory. If she uses the global category test, this is a straightforward determination.

The global category test cannot be manipulated, in the sense that every random device $\zeta$ fails this test with certainty on all paths except a first category set of them.

The paper is organized as follows. In Section 2, we introduce our basic concepts and show some classic examples of tests (calibration and likelihood tests). In Section 3, we construct the global category test and provide an informal discussion of nonmanipulability. In Section 4, we show that the global category test can be modified to address the problem that, in practice, testers observe only finite data sets. We also show that the global category test can be modified so as to belong to the familiar class of likelihood tests. In addition, we exhibit a large class of tests that are not necessarily prequential tests, but can be manipulated. Finally, still in Section 4, we offer some additional results on the manipulability of random prequential tests and an informal discussion on the implications (for the prequential principle) of the finding that prequential tests can be manipulated while some nonprequential tests cannot be manipulated. Proofs are relegated to the Appendix.

**2. Basic concepts.** Each period one outcome, 0 or 1, is observed (our results generalize to any finite number of outcomes per period). Let $\Omega = \{0,1\}^\infty$ be the set of all *paths*, that is, infinite histories. A finite history $s_m \in \{0,1\}^m$, $m \geq t$, (or a path $s \in \Omega$) is an *extension* of a history $s_t \in \{0,1\}^t$ if the first $t$ outcomes of $s_m$ or $s$ coincide with the outcomes of $s_t$. In the opposite direction, let $s_m \mid t$ (or $s \mid t$) be the history $s_t \in \{0,1\}^t$ whose outcomes coincide with the first $t$ outcomes of $s_m$ or $s$. A *cylinder* with base on $s_t$ is the set $C(s_t) \subset \{0,1\}^\infty$ of all infinite extensions of $s_t$. We endow $\Omega$ with the topology that compares unions of cylinders with a finite base. Let $\Im_t$ be the algebra that consists of all finite unions of cylinders with base on $\{0,1\}^t$. Denote by $N$ the set of natural numbers. Let $\Im$ be the $\sigma$-algebra generated by the algebra $\Im^0 := \bigcup_{t \in N} \Im_t$; that is, $\Im$ is the smallest $\sigma$-algebra that contains $\Im^0$.

Let $\Delta(\Omega)$ be the set of all probability measures on $(\Omega, \Im)$. We endow $\Delta(\Omega)$ with the weak*-topology and the $\sigma$-algebra of Borel sets (i.e., the smallest $\sigma$-algebra that contains all open sets in weak*-topology).

As is well known, the weak*-topology consists of all unions of finite intersections of sets of the form

$$\{Q \in \Delta(\Omega) : |E^P h - E^Q h| < \varepsilon\},$$

where $E$ stands for the expected-value operator $P \in \Delta(\Omega)$, $\varepsilon > 0$, and $h$ is a real-valued and continuous function on $\Omega$. We refer the reader to Rudin (1973) for additional details on the weak*-topology. We also define $\Delta\Delta(\Omega)$ as the set of probability measures on $\Delta(\Omega)$.

Before any data are observed, an expert named Bob announces a probability measure $P \in \Delta(\Omega)$, which (Bob claims) describes how nature will generate the data. To simplify the language, we call a probability measure a *theory*. A tester named Alice tests Bob's theory empirically.

DEFINITION 1.   A test is a function $T : \Omega \times \Delta(\Omega) \to \{0, 1\}$.

That is, a test is defined as an arbitrary function that takes as input a theory and a path and returns a verdict that is 0 or 1. When the test returns a 1, it does not reject (or, simply, accepts) the theory. When a 0 is returned, the theory is rejected.

Any test divides paths into those in $A_P := \{s \in \Omega \mid T(s, P) = 1\}$, where the theory $P$ is accepted, and those in $\Omega - A_P$, where the theory is rejected. The set $A_P$ is called the *acceptance set*, and its complement $\Omega - A_P$ is called the *rejection set.* We consider only tests $T$ such that the acceptance sets $A_P$ are $\Im$-measurable.

Given a theory $P \in \Delta(\Omega)$, a path $s \in \Omega$ and a history $s_t = s \mid t$, let

$$f_0^P(s) := P(C(1)) \quad \text{and} \quad f_t^P(s) := \frac{P(C(s_t, 1))}{P(C(s_t))}$$

be forecasts made along $s$. $f_t^P(s)$ is arbitrarily defined as 0 when $P(C(s_t)) = 0$. The forecasts of $P$ and $P'$ are *equivalent along* $s$ if $f_t^P(s) = f_t^{P'}(s)$ for all periods $t \geq 0$. A test $T$ is *prequential* if for any given two theories $P$ and $P'$, equivalent along $s$, $T(s, P) = 0$ if and only if $T(s, P') = 0$. So, a prequential test rejects or accepts a theory based only on the forecasts made by the theory along the realized path. Fix any $\varepsilon \in [0, 1]$.

DEFINITION 2.   A test $T$ does not reject the data-generating process with probability $1 - \varepsilon$ if for any $P \in \Delta(\Omega)$

$$P(A_P) \geq 1 - \varepsilon.$$

That is, a test does not reject the data-generating process (with high probability) if, no matter which probability measure $P$ generates the data, $P$ is not likely to be rejected according to its own probability distribution over paths.

2.1. *Classic examples of tests.* Given a path $s \in \Omega$, let $I_t(s)$ be the $t$th outcome of $s$. The test

$$T(s, P) = 1 \quad \text{if and only if} \quad \lim_{n \to \infty} \frac{1}{n} \sum_{t=1}^{n} [f_{t-1}^{P}(s) - I_t(s)] = 0$$

requires the average forecast of 1 to match the empirical frequency of 1.

The test in Lehrer ([2001](#)) requires the match between average forecasts and empirical frequencies to occur on several subsequences. The calibration test requires the empirical frequency of 1 to be near $p \in [0, 1]$ in the periods in which 1 was forecasted with probability close to $p$.

One well-known class of tests not based on matching empirical frequencies is likelihood tests. Let $Q : \Delta(\Omega) \to \Delta(\Omega)$ be any function that takes a theory $P$ as an input and returns an alternative theory $Q_P$ as an output. The *likelihood test* $T$ is defined by

$$T(s, P) = 1 \qquad \text{iff } \forall_{n \in N} P(C(s_n)) \neq 0$$

and

$$\limsup_{n \to \infty} \frac{Q_P(C(s_n))}{P(C(s_n))} < \infty, \qquad s_n = s \mid n.$$

This test requires that the likelihood of $Q_P$ does not become arbitrarily larger than the likelihood of $P$. A proof that calibration and likelihood tests do not reject the data-generating process with probability one can be found in Dawid ([1982](#), [1985](#)).

Calibration tests are prequential tests because, like many other well-known tests, they take as an input not the entire theory but only the data and the forecasts made by the theory along the realized sequence of outcomes. A likelihood test may or may not be a prequential test, depending upon the way the alternative theory $Q_P$ is selected as a function of $P$. In Section [4](#), we exhibit examples of prequential and nonprequential likelihood tests.

2.2. *Manipulating tests.* Foster and Vohra ([1998](#)) show that Bob can pass the calibration test on all paths if he is allowed to select theories at random. More precisely, the calibration test can be manipulated with probability one according to the following definition:

DEFINITION 3. A test $T$ can be manipulated with probability $1 - \varepsilon$ if there exists a random generator of theories $\zeta_T \in \Delta\Delta(\Omega)$ such that, for every path $s \in \Omega$,

$$\zeta_T(\{P \in \Delta(\Omega) \mid T(s, P) = 1\}) \geq 1 - \varepsilon.$$

Naturally, our definition of manipulability requires a measurability provision on the sets $\{P \in \Delta(\Omega) : T(s, P) = 1\}$. However, to ease the exposition, we deal with measurability issues in the Appendix.

Fudenberg and Levine (1999), Lehrer (2001), Sandroni, Smorodinsky and Vohra (2003) show that generalized calibration tests can also be manipulated. Sandroni (2003), Vovk and Shafer (2005), Olszewski and Sandroni (2008) and Shmaya (2008) show that prequential tests that can be manipulated. A partial review of this literature can be found in Cesa-Bianchi and Lugosi (2006). We also refer the reader to Gneiting, Balabdaoui and Raftery (2007) for additional comments on this literature.

The random generator $\zeta_T$ may depend only on the test $T$. It follows that, even if Bob does not know the data-generating process, he can be very confident that by selecting theories at random before any data are observed he will not fail any of these tests, no matter which data are realized in the future.

The crucial property we seek is the existence, for every random generator of theories, of a path at which rejection may occur (because it ensures that it is feasible to reject strategically produced theories). However, stronger properties may be obtained. Alice may be interested in paths such that randomly produced theories fail the test with near certainty (as opposed to probability higher than $\varepsilon$) and such that the sets of these paths are larger than a single path.

DEFINITION 4. Fix a test $T$. Given a random generator of theories $\zeta \in \Delta\Delta(\Omega)$ and $\varepsilon \geq 0$, let $R_\zeta^\varepsilon \subseteq \Omega$ be the set of all paths $s \in \Omega$ such that

$$\zeta(\{P \in \Delta(\Omega) \mid T(s, P) = 0\}) \geq 1 - \varepsilon.$$

The set $R_\zeta^\varepsilon$ is called the *revelation set*, where the random generator of theories $\zeta$ fails to manipulate the test with probability $1 - \varepsilon$.

In Section 3, we exhibit a test such that for all random generators of theories $\zeta$, the revelation sets are always non-empty (and topologically large).

**3. A nonmanipulable test.** In this section, we construct a test that cannot be manipulated, that is, a test with nonempty revelation sets. In addition, we show that the revelation sets are topologically large in the following sense: Given a subset $A$ of a (complete metric) space, let $\bar{A}$ be the closure of $A$. A set $A$ is called *nowhere dense* when the interior of its closure $\bar{A}$ is empty. A *first-category* set is a countable union of nowhere-dense sets. A first-category set may be regarded as (topologically) small. The complement of a first-category set may be regarded as (topologically) large. We refer the reader to Oxtoby (1980) for these definitions and some basic results regarding first-category sets.

We now define our test. Let $S = \{s^1, s^2, \ldots\}$ be a countable dense subset of $\Omega$ (i.e., for every cylinder $C$ with finite base, there exists a path $s^i$ that belongs to $C$; e.g., $C$ can be the set of paths with all but a finite number of outcomes equal to zero). Fix any $k \in N$. For every path $s^i$, there exists a period $t \in N$ such that the cylinder $C(s_t^i)$ with base on the finite history $s_t^i = s^i \mid t$ satisfies

$$(3.1) \qquad\qquad P(C(s_t^i) - \{s^i\}) \leq \frac{1}{2^{k+i}}.$$

Indeed, the sequence of sets $C(s_t^i) - \{s^i\}$ is descending (as $t$ goes to infinity), and its intersection is empty. So, $P(C(s_t^i) - \{s^i\})$ goes to zero as $t$ goes to infinity. Let $t(i, k, P)$ be the smallest natural number such that (3.1) is satisfied.

Given a dense set of paths $S = \{s^1, s^2, \ldots\}$, the *global category test* $\widehat{T}_S$ can be defined as follows:

$$\widehat{T}_S(s, P) = \begin{cases} 0, & \text{if } s \notin S \text{ and } s \in \bigcap_{k=1}^{\infty} \bigcup_{i=1}^{\infty} C(s_{t(i,k,P)}^i); \\ 1, & \text{otherwise.} \end{cases}$$

Let $\widehat{R}_\zeta^0$ be the revelation set for the random generator of theories $\zeta \in \Delta\Delta(\Omega)$.

THEOREM 1. *Fix any countable dense set $S \subset \Omega$. The global category test $\widehat{T}_S$ does not reject the data-generating process with probability one. Given any random generator of theories $\zeta \in \Delta\Delta(\Omega)$, the revelation set $\widehat{R}_\zeta^0$ is the complement of a first-category set of paths.*

The global category test controls the type I error of rejecting the data-generating process and it cannot be manipulated. If the data-generating process is announced, then it passes the global category test with probability one. On the other hand, no matter which random generator of theories $\zeta$ is employed, failure is inevitable on the paths of the revelation set $\widehat{R}_\zeta^0$. This set is nonempty (and topologically large).

The main objective of this paper is to construct a test that is nonmanipulable according to Definition 3, that is, in the sense studied in the existing literature. Naturally, our result is limited to our definition of manipulation. Alternative definitions of strategic manipulation are beyond the scope of this paper.

The global category test can be combined with any other test to produce a harder test. Even so, no test can avoid the following difficulty: for any finite collection of probability measures $Q_1, \ldots, Q_k$ over outcome paths, and for any test $T$ that does not reject the data-generating process with high probability, Bob can ensure that his randomly selected theory is unlikely to be rejected on a set of outcome paths whose $Q_j$-probabilities ($j = 1, \ldots, k$) are arbitrarily close to one [see Sandroni and Olszewski (2008)].

Since the proof of Theorem 1 is somewhat involved, it is relegated to the Appendix, but the key idea of this proof is simple. The rejection sets for the global category test have been defined as an intersection of unions of small cylinders $C_1$, $C_2, \ldots$ around paths $s^1$, $s^2, \ldots$ (from which the paths $s^1$, $s^2, \ldots$ have additionally been removed). The fact that a cylinder around each path is included in the rejection set guarantees that, for any single theory, the rejection set is topologically large; moreover, the fact that those cylinders are small guarantees that the data-generating process will not be rejected.

Suppose now that we are given a random generator of theories $\zeta$. Since cylinders $C_1$, $C_2, \ldots$ are defined for every single theory around the fixed set of paths $s^1$, $s^2, \ldots$, it follows that sufficiently small cylinders around those paths $s^1$, $s^2, \ldots$ will be contained in cylinders $C_1$, $C_2, \ldots$, respectively, for a set of theories whose $\zeta$-probability is sufficiently close to one. As the intersection of topologically large sets is itself a topologically large set, there exists a topologically large set contained in the rejection sets for a set of theories whose $\zeta$-probability is equal to one.

3.1. *Informal discussion on manipulability.* A result concerning the manipulability of a test can be interpreted in different ways. Consider a forecaster whose objective is to be calibrated. The Foster and Vohra (1998) result shows that the forecaster's goal can be achieved no matter how the data evolves in the future. Hence, from this forecaster's perspective, the Foster and Vohra (1998) result is positive. The results of Sandroni (2003), Vovk and Shafer (2005), Olszewski and Sandroni (2008) and Shmaya (2008) show that, like calibration, many observable properties of the data-generating process can be obtained by a strategic forecaster, no matter how the data evolves in the future. So, these results show how to produce forecasts that, in the future, will prove to have some observable properties of the data-generating process. However, these strategically produced forecasts may remain bounded away from the predictions of the data-generating process, and these forecasts need not have all observable properties of the data-generating process simultaneously (only the properties used to define some specific manipulable test). Hence, whether strategically produced forecasts are desirable from the perspective of a forecaster is an arguable point, as this may depend upon the objective of the forecaster. Now, consider the concept of manipulability from the viewpoint of a tester.

Assume that Alice wants to determine whether Bob has information about the data-generating process that she does not have. If this is Alice's objective, then a manipulable test (e.g., the calibration test) has limited use to her when Bob knows the test that will be used (see Section 4.2 for a discussion of the case in which Bob does not know the test at the time of delivering his theory). Even in the extreme case that Bob is completely uninformed about the data-generating process, he can strategically pass such a test. Hence, Alice knows from the outset (i.e., before any data is revealed) the verdict that, with near certainty, the test will deliver once the data is revealed.

The difficulty with manipulable tests can be understood in the context of a contracting problem between a tester (Alice) and an expert (Bob) who claims, before any data is observed, to know the data-generating process. Alice does not know the data-generating process, and she is willing to pay a reward that gives Bob utility $u > 0$ if he announces his theory at period zero. In order to discourage Bob from delivering an arbitrary theory, Alice stipulates a penalty if Bob's theory is rejected by her test. This penalty gives Bob disutility $d > u$ (Bob does not discount the future). Bob observes Alice's test before deciding whether to accept Alice's contract. Bob receives no reward and no penalty if he does not accept the contract.

Alice looks for a screening contract, which will be accepted by Bob if informed about the data-generating process and rejected by Bob if uninformed about the data-generating process. Then, Alice learns whether Bob is informed from Bob's choice on whether to accept the contract.

Assume that the test does not reject the data-generating process with probability $1 - \varepsilon$, where $\varepsilon$ is small enough so that $u - d\varepsilon > 0$ and $u - d(1 - \varepsilon) < 0$. Bob, if informed, accepts the contract because his expected utility with the contract is positive (and without the contract his utility is zero). On the other hand, if Bob is uninformed, then he faces uncertainty. He does not know the odds that any theory is rejected. Assume that Bob evaluates his prospects based on the minimum expected utility he obtains. This is the most pessimistic behavioral rule of decision under uncertainty among those axiomatized by Gilboa and Schmeidler ([1989](#)). Formally, if Bob is uninformed and uses a random generator of theories $\zeta \in \Delta\Delta(\Omega)$, then his payoff, with a contract, is

$$(3.2) \qquad u - d \sup_{s \in S} \zeta \{ P \in \Delta(\Omega) \mid T(s, P) = 0 \}.$$

If Alice's test is manipulable, then Bob, although extremely averse to uncertainty, accepts her contact because his payoff with the contract (3.2) is positive. Hence, a manipulable tests cannot be used to construct a screening contract that Bob, if informed, accepts and if uninformed, does not accept.

On the other hand, if Alice's test is nonmanipulable, then Bob, when uninformed, does not accept Alice's contract, as his payoff with the contract is negative. So, a nonmanipulable test does produce a screening contract. Bob accepts this contract (and, therefore, Alice pays Bob) only when Bob delivers to Alice something she values (the data-generating process), as opposed to a theory selected by a method that she can produce on her own without having to pay an expert for it.

3.2. *Limitations of the main result.*   The main objective of this paper is
to construct a nonmanipulable test. Nonmanipulability is in our opinion a
desirable property (at least as far as Bob knows in advance the test that will
be used), because it ensures the feasibility of the rejection of a strategic, but
otherwise completely uninformed, expert. However, we do not claim that an
uninformed expert is likely to be rejected by a global category test, nor do
we claim that a theory that passes our test must be closely related to the
data-generating process.

In addition, we are not confidently advocating the use of global category
tests. An uninformed expert may accept Alice's contract, even if she adopts
a global category test, when Bob adopts a less pessimistic behavioral rule
to evaluate uncertain prospects. Moreover, global category tests have not
been extensively studied, and they may have undesirable properties. On the
other hand, global category tests do not reject the data-generating process
with probability one. Hence, they can be combined with any other test
without reducing the odds that the data-generating process is rejected. This
combined test remains nonmanipulable.

**4. Additional results.**   In this section, we show modified forms of the
global category test and discuss some of its properties. We begin with the
simple observation that tests, as defined in Section 2, give a verdict only after
observing an infinite history. In practice, a tester can only observe finite data
sets and, therefore, may find more applicable tests that give some verdict
with a finite number of outcomes.

DEFINITION 5.   Rejection tests have the property that, for any theory
$P \in \Delta(\Omega)$, the rejection set $\Omega - A_P$ is a union of cylinders.

So, a rejection test rejects a theory in finite time.

DEFINITION 6.   A test $T_2$ is harder than a test $T_1$ if

$$\{s \in \Omega : T_2(s, P) = 1\} \subseteq \{s \in \Omega : T_1(s, P) = 1\}.$$

If $T_2$ is harder than $T_1$, then rejection by $T_1$ implies rejection by $T_2$.

PROPOSITION 1.   *Fix any* $\delta \in (0, 1 - \varepsilon]$. *Let* $T_1$ *be any test that does
not reject the data-generating process with probability* $1 - \varepsilon$. *There exists
a rejection test* $T_2$ *that is harder than the test* $T_1$ *and does not reject the
data-generating process with probability* $1 - \varepsilon - \delta$.

Proposition 1 is a direct corollary of the following well-known result: for
any given probability measure $P \in \Delta(\Omega)$ and $\delta > 0$, any set $A \in \Im$ can be enlarged to an open set $U \supset A$ such that $P(U) < P(A) + \delta$ [see Ulam's theorem, Theorem 7.1.4 in Dudley (1989)].

Proposition 1 shows that, given any test $T_1$, there exists a rejection test $T_2$ that does not reject the data-generating process with probability almost as high as $T_1$, and the acceptance sets of $T_2$ are contained in those of $T_1$. Thus, $T_2$ is more capable than $T_1$ to reject theories. For example, if $T_1$ has nonempty (or topologically large) revelation sets, then $T_2$ is a rejection test with nonempty (or topologically large) revelation sets. We refer to $T_2$ as the rejection test associated with the test $T_1$.

Given $\varepsilon > 0$, let $\hat{T}^\varepsilon$ be a rejection test (that does not reject the data-generating process with probability $1 - \varepsilon$) associated with a global category test $\hat{T}_S$. Corollary 1 follows immediately from Proposition 1 and Theorem 1.

COROLLARY 1. *The rejection test $\hat{T}^\varepsilon$ does not reject the data-generating process with probability $1 - \varepsilon$, and it cannot be manipulated. The revelation sets of $\hat{T}^\varepsilon$ contain the revelation sets of the global category test, and, hence, their complements are first-category sets.*

For any given theory $P$, the test $\hat{T}^\varepsilon$ delivers the rejection set of the theory $P$ that consists of finite histories of outcomes that must be regarded as sufficiently inconsistent with the theory $P$ to justify its rejection. This test does not reject the actual data-generating process and, in contrast with manipulable tests, maintains the possibility of theory rejection (in finite time) even if an uninformed individual randomizes with the intent of manipulating the test. In addition, it follows directly from our proof that each revelation set of $\hat{T}^\varepsilon$ contains an open and dense set. Hence, the revelation sets of $\hat{T}^\varepsilon$ contain finite-histories.

4.1. *A nonmanipulable likelihood test.* In this section, we modify a global category test so that, in this modified form, it belongs to the familiar class of likelihood tests (and remains nonmanipulable).

Let an *atom* of a theory $P$ be a path that $P$ assigns strictly positive measure. Let $\mathcal{A} \subseteq \Delta(\Omega)$ be the set of all theories that attach probability one to a set comprising finitely many atoms.

Let $\hat{T}_S$ be a global category test. Given any natural number $m \geq 1$, let $\hat{T}^m$ be a rejection test that is harder than $\hat{T}_S$ and does not reject the data-generating process with probability $\frac{1}{2(m+1)^3}$. Given a theory $P$, let $\hat{R}_P^m$ be the rejection set of $P$ (for test $\hat{T}^m$). So, by definition,

$$P(\hat{R}_P^m) \leq \frac{1}{2(m+1)^3}.$$

Given any theory $P \in (\mathcal{A})^c$, let $C_P^m$ be a cylinder with finite base such that

$$(4.1) \qquad\qquad 0 < P(C_P^m) < \frac{1}{2(m+1)^3}.$$

The existence of a cylinder that satisfies (4.1) is shown in Remark A.1 in the Appendix. Let $R_P^m := \hat{R}_P^m \cup C_P^m$, and let $Q_P^m$ be the theory $P \in (\mathcal{A})^c$ conditional of $R_P^m$; that is, for every $A \in \Im$,

$$(4.2) \qquad\qquad Q_P^m(A) = \frac{P(A \cap R_P^m)}{P(R_P^m)}.$$

Let $\pi(m) = \frac{1}{(m+1)m}$. Let $Q_P$ be the theory defined by

$$Q_P = \sum_{m=1}^{\infty} \pi(m) Q_P^m.$$

Given a theory $P \in (\mathcal{A})^c$, $Q_P$ is well defined because

$$\sum_{m=1}^{\infty} \pi(m) = 1.$$

The function $Q \colon \Delta(\Omega) \longrightarrow \Delta(\Omega)$, such that $Q(P) = Q_P$ if $P \in (\mathcal{A})^c$ and $Q(P) = P$ if $P \in \mathcal{A}$, defines the likelihood test $\bar{T}$.

PROPOSITION 2.  *The likelihood test $\bar{T}$ is harder than the global category test $\hat{T}_S$.*

The likelihood test $\bar{T}$, like any other likelihood test, does not reject the data-generating process with probability 1. It follows from Proposition 2 that the revelation sets of the likelihood test $\bar{T}$ contain those of the global category test $\hat{T}_S$. Hence, the likelihood test $\bar{T}$ is nonmanipulable.

4.2. *The prequential principle.*  Prequential tests are manipulable, but some nonprequential tests (e.g., the global category tests) are nonmanipulable. These results pose a difficulty for the prequential principle, because they indicate that the prequential principle must be discarded to produce a nonmanipulable test. We examine in closer detail the relationship between manipulability and the prequential principle. We now consider random tests. That is, we allow Alice to randomize among prequential tests that do not reject the data-generating process.

Let $(\Theta, \mathcal{B}, \tilde{v})$ be a probability space where $\Theta$ is a parameter space, $\mathcal{B}$ is a $\sigma$-algebra and $\tilde{v}$ is a probability measure on $(\Theta, \mathcal{B})$. A *random test* is probability space $(\Theta, \mathcal{B}, \tilde{v})$ and a function $\tilde{T} \colon \Theta \times \Omega \times \Delta(\Omega) \to \{0, 1\}$. So, for every parameter $\theta \in \Theta$, $\tilde{T}_\theta : \Omega \times \Delta(\Omega) \to \{0, 1\}$, $\tilde{T}_\theta(s, P) = \tilde{T}(\theta, s, P)$ is a test. We also define $\tilde{T}_s(\theta, P) = \tilde{T}(\theta, s, P)$ and $\tilde{T}_P(\theta, s) = \tilde{T}(\theta, s, P)$. We assume that $\tilde{T}$ is measurable, jointly, with respect to $\theta$, $s$ and $P$ (and we refer to this assumption as the joint measurability condition).

A random test is prequential if, for every $\theta \in \Theta$, $\tilde{T}_\theta$ is a prequential test. A random test does not reject the data-generating process with probability one if, for every $P \in \Delta(\Omega)$,

$$(\tilde{v} \times P)(\{(\theta, s) \in \Theta \times \Omega : \tilde{T}_P(\theta, s) = 1\}) = 1.$$

We now construct a prequential random test $\tilde{T}$ that does not reject the data-generating process with probability one. Given a countable dense set $S = \{s^1, s^2, \ldots\}$, let $Q^S \in \Delta(\Omega)$ be a probability distribution defined by

$$Q^S(s^i) = \frac{1}{i(i+1)}$$

for $i = 1, 2, \ldots$ and $Q^S(\Omega - S) = 0$. That is, $Q^S$ assigns full measure to $S$. Since

$$\sum_{i=1}^{\infty} \frac{1}{i(i+1)} = 1,$$

$Q^S$ is well defined.

Let $T_S^{LR}$ be the likelihood test such that, for every theory $P$, the alternative theory is $Q^S$. This test is prequential because the alternative theory is fixed independently of the theory announced by Bob. By general properties of likelihood ratio tests, $T_S^{LR}$ does not reject the data-generating process with probability 1 [see Dawid (1982)].

Let $\Theta = \Omega = \{0, 1\}^\infty$ be the parameter space. Given a path $\theta \in \Theta$, let $S_\theta \subset \Omega$ be the set of all paths that coincide with $\theta$ in all but a finite number of periods. That is,

$$S_\theta = \{s \in \Omega : |\{t : I_t(s) \neq I_t(\theta)\}| < \infty\}$$

$$= \bigcup_{t=1}^{\infty} \{0, 1\}^{t-1} \times \{I_t(\theta)\} \times \{I_{t+1}(\theta)\} \times \cdots ;$$

of course the set $S_\theta$ is countable and dense in $\Omega$. We define the *randomized likelihood test* $\tilde{T}^{LR}$ as follows. First, let $\dot{\Delta}(\Omega) \subseteq \Delta(\Omega)$ be the set of all theories that assigns zero measure to any single path. Alice draws randomly a path $\theta \in \Omega$ according to a probability measure $\tilde{v} \in \dot{\Delta}(\Omega)$ and then tests Bob with the likelihood test $T_{S_\theta}^{LR}$. So, $\tilde{T}^{LR}(\theta, s, P) = T_{S_\theta}^{LR}(s, P)$.

To make this definition precise, we must say how we represent the set $S_\theta$ as a sequence $s_\theta^1, s_\theta^2, \ldots$, because the probability distribution $Q^{S_\theta}$ depends on the order of the paths in $S_\theta$. In addition, if we took an arbitrary representation, the test $\widetilde{T}^{LR}$ could potentially violate our joint measurability condition.

One possible way to represent the set $S_\theta$ as a sequence of paths $s_\theta^1, s_\theta^2, \ldots$ is such that for every $t = 1, 2, \ldots$ the paths from $\{0, 1\}^t \times \{I_t(\theta)\} \times \{I_{t+1}(\theta)\} \times \cdots$ precede the paths from $\{0, 1\}^t \times \{I_{t+1}(\theta)\} \times \{I_{t+2}(\theta)\} \times \cdots - \{0, 1\}^{t-1} \times \{I_t(\theta)\} \times \{I_{t+1}(\theta)\} \times \cdots$, and paths from $\{0, 1\}^t \times \{I_{t+1}(\theta)\} \times \{I_{t+2}(\theta)\} \times \cdots - \{0, 1\}^{t-1} \times \{I_t(\theta)\} \times \{I_{t+1}(\theta)\} \times \cdots$ are ordered lexicographically. As shown in the Appendix (Lemma A.2), the joint measurability condition is satisfied with the sets $S_\theta$ ordered in this way.

Bob knows that he is tested according to this protocol, but he does not know the test selected by Alice. Given any random generator of theories $\zeta \in \Delta\Delta(\Omega)$, let $R_\zeta^0(\theta)$ be the revelation set of the test $T_{S_\theta}^{LR}$.

PROPOSITION 3. *Consider a randomized likelihood test $\widetilde{T}^{LR}$, $(\Omega, \Im, \tilde{v})$, $\tilde{v} \in \dot{\Delta}(\Omega)$. For every random generator of theories $\zeta \in \Delta\Delta(\Omega)$, $\tilde{v}$-almost surely, $R_\zeta^0(\theta)$ is the complement of a first-category set of paths.*

Assume that Alice tests Bob with a randomized likelihood test and that Bob uses an arbitrary random generator of theories $\zeta$. Proposition 3 shows that, with $\tilde{v}$-probability one, Alice selects a prequential test $T_{S_\theta}^{LR}$ for which there *exist* a topologically large set of paths that, if realized, $\zeta$-almost surely, reject Bob's theory.

We now provide a general result showing that random prequential tests, including the randomized likelihood test, are manipulable as far as a natural generalization (to random tests) of definition 3 of manipulability goes.

PROPOSITION 4. *Fix any $\delta > 0$. Let $(\Theta, \mathcal{B}, \tilde{v})$ and $\widetilde{T}$ be a prequential random test (satisfying the joint measurability condition) that does not reject the data-generating process with probability 1. Then, there exists a random generator of theories $\tilde{\zeta}$ such that, on any path $s \in \Omega$,*

$$\tilde{\zeta}(\{P \in \Delta(\Omega) : \tilde{v}\text{-almost surely } \widetilde{T}(\theta, s, P) = 1\}) \geq 1 - \delta.$$

The random generator of theories $\tilde{\zeta}$ may depend upon $\delta$, $(\Theta, \mathcal{B}, \tilde{v})$ and $\widetilde{T}$, but it does not require (for its construction) any distributional assumptions over the future realizations of the paths. Proposition 4 shows that Bob can produce theories according to a random device $\tilde{\zeta}$ such that, no matter which path $s$ is realized, it is unlikely (odds given by $\tilde{v}$ and $\tilde{\zeta}$) that Alice selects a test that rejects Bob's theory.

To reconcile Propositions 3 and 4, consider the randomized likelihood test $(\Omega, \Im, \tilde{v})$, $\widetilde{T}^{LR}$. This test satisfies the joint measurability conditions, and so, the conclusions of Propositions 3 and 4 hold for this test. Let $E^{\tilde{v}}$ and $E^\zeta$ be the expectation operators associated with $\tilde{v}$ and $\zeta$, respectively.

By Proposition 3, for every $\zeta \in \Delta\Delta(\Omega)$,

$$(4.3) \qquad E^{\tilde{v}}\left\{\inf_{s \in \Omega} E^{\zeta}\{\tilde{T}_s^{LR}\}\right\} = 0.$$

By Proposition 4, for every $\delta > 0$ there exists $\tilde{\zeta} \in \Delta\Delta(\Omega)$ such that

$$(4.4) \qquad \inf_{s \in \Omega} E^{\tilde{v}} E^{\tilde{\zeta}}\{\tilde{T}_s^{LR}\} \geq 1 - \delta.$$

From (4.3), Bob knows that Alice (almost surely) selects a test for which there are paths that, if realized, will (almost surely) reject Bob's randomly selected theories. From (4.4), Bob knows that if he selects a theory with carefully designed odds (by $\tilde{\zeta}$), then (no matter how the data evolves in the future) it is unlikely that the selected theory will be rejected by Alice's test.

Consider the question of whether a strategic, but uninformed, expert can pass a random prequential test (which does not reject the data-generating process). Proposition 3 seems to answers this question in the negative while Proposition 4 seems to answer this question in the positive. Hence, results (4.3) and (4.4) leave room for different interpretations. In our viewpoint, Proposition 4 is a natural generalization (to random tests) of existing results showing that prequential tests are manipulable. Our preference (for result 4 over 3) can be understood in the context of decision-making under uncertainty as described in Section 3.1.

Assume that Alice offers a contract to Bob. If Bob accepts the contract, he delivers a theory to Alice and receives positive payment. However, if Bob's theory is rejected by the test selected by the randomized likelihood test, then Bob is penalized. Now, assume that Bob knows nothing about the data-generating process. Then, Bob faces uncertainty about the probabilities of the future realizations of the data. By definition, Bob knows the odds that Alice uses to select her test, and he also knows the odds that he uses to select his theory. Hence, with regards to theory selection and to test selection, Bob faces common risk. The most pessimistic behavioral rule of decision under uncertainty, among those axiomatized by Gilboa and Schmeidler (1989), determines Bob's prospects by his expected utility computed in the worse-case scenario [as in (4.4)]. By this rule of decision under uncertainty, Bob accepts Alice's contract. Hence, prequential random tests do not screen informed and uniformed experts.

The results showing that prequential tests (and even randomizations over prequential tests, Proposition 3 withstanding) are manipulable, combined with the fact that some nonprequential tests are nonmanipulable, poses a difficulty for the prequential principle (as far as testing potentially strategic experts goes). However, this difficulty may not persist under conditions that are beyond the scope of this paper. Fortnow and Vohra (2007) show a prequential test that is computationally demanding to manipulate. In addition, Olszewski and Sandroni (2008) show a prequential test that cannot be manipulated when the domain of permissible theories (i.e., the theories the expert is allowed to announce) is restricted. It is not known whether the results in this paper extend to the case of multiple experts [see Al-Najjar and Weinstein (2007) and Feinberg and Stewart (2007) for some results on testing several experts simultaneously]. Finally, while any given prequential test (which does not reject the data-generating process) can be manipulated, it is not possible to manipulate all prequential tests (that do not reject the data-generating process) simultaneously [see Olszewski and Sandroni (2008)].

4.3. *Nonprequential manipulable tests.* As mentioned in the Introduction, prequential tests are manipulable. So, prequentiality is a sufficient condition for manipulability, but it is not a necessary condition. We now show several tests, many of them nonprequential, that can be manipulated.

DEFINITION 7. Acceptance tests have the property that, for any theory $P \in \Delta(\Omega)$, the acceptance set $A_P$ is a union of cylinders.

In an acceptance test, the acceptance sets are open. In a rejection test, the rejection sets are open.

PROPOSITION 5. *Fix any $\varepsilon \in [0,1]$ and $\delta \in (0, 1 - \varepsilon]$. Let $T$ be an acceptance test that does not reject the data-generating process with probability $1 - \varepsilon$. Then, the test $T$ can be manipulated with probability $1 - \varepsilon - \delta$.*

A formal proof of Proposition 5 is presented in the Appendix. An intuition is as follows: let $V : \Delta(\Omega) \times \Delta\Delta(\Omega) \to [0,1]$ be a function defined by $V(P, \zeta) = E^P E^\zeta T$; that is, $V(P, \zeta)$ is the probability of the verdict 1 if $P$ is the data-generating process and $\zeta$ is the random generator of theories used by Bob. By assumption, for every $P \in \Delta(\Omega)$, there exists $\zeta_P \in \Delta\Delta(\Omega)$ (a deterministic generator of theories that assigns probability one to $P$) such that $V(P, \zeta_P) = 1 - \varepsilon$. Thus, if the conditions of Fan's minmax theorem are satisfied, then there also exists $\zeta_T \in \Delta\Delta(\Omega)$ such that $V(P, \zeta_T) \geq 1 - \varepsilon - \delta$ for every $P \in \Delta(\Omega)$. This yields the result, since $V(P, \zeta_T) = \zeta_T(\{Q \in \Delta(\Omega) \mid T(s, Q) = 1\})$ if $P$ is the degenerated measure that assigns probability one to $s$.

As is well known, $\Delta(\Omega)$ is compact in the weak*-topology and $V$ is a bilinear function. Hence, the conditions of Fan's minmax theorem (see the Appendix for this result) are satisfied if $V$ is lower semi-continuous with respect to $P$. We show that this lower semi-continuity follows from the openness of acceptance sets. It is here that the assumption that $T$ is open turns out to be essential. In the case of a rejection test, $V$ is upper semi-continuous with respect to $P$, but not necessarily lower semi-continuous.

## APPENDIX: PROOFS

PROOF OF THEOREM 1. Define $A_k^i(P) := \Omega - C(s_{t(i,k,P)}^i)$ as the complement of $C(s_{t(i,k,P)}^i)$. Let

$$A_k(P) := \bigcap_{i=1}^{\infty} A_k^i(P) \quad \text{and} \quad \hat{A}_P := \bigcup_{k=1}^{\infty} A_k(P) \cup \bigcup_{i=1}^{\infty} \{s^i\}.$$

Each set $A_k(P)$ is an intersection of closed sets and is therefore closed itself. By construction, each set $A_k(P)$ has an empty interior; indeed, the complement of $A_k(P)$; that is, the set

$$\Omega - A_k(P) = \bigcup_{i=1}^{\infty} C(s_{t(i,k,P)}^i)$$

is open and dense. Hence, $\hat{A}_P$ is a first-category set. Notice that $\hat{A}_P$ is the acceptance set of $P \in \Delta(\Omega)$.

Since

$$\Omega - \hat{A}_P \subset \Omega - \left( A_k(P) \cup \bigcup_{i=1}^{\infty} \{s^i\} \right) \subset \bigcup_{i=1}^{\infty} [C(s_{t(i,k,P)}^i) - \{s^i\}]$$

for all $k \in N$,

$$P(\Omega - \hat{A}_P) \leq \sum_{i=1}^{\infty} P(C(s_{t(i,k,P)}^i) - \{s^i\}) \leq \sum_{i=1}^{\infty} \frac{1}{2^{k+i}} = \frac{1}{2^k},$$

also for all $k \in N$, which yields that $P(\hat{A}_P) = 1$. Thus, the global category test does not reject the data-generating process with probability one. It remains to show that the test cannot be manipulated. Suppose we are given a $\zeta \in \Delta\Delta(\Omega)$. We first show that there exists a subset $\hat{A}$ of $\Omega$, which is a countable union of closed sets with empty interior, and a Borel set $B \subset \Delta(\Omega)$ such that

(A.1) $$\zeta(B) = 1$$

and

(A.2) $$\forall_{P \in B} \qquad \hat{A}_P \subset \hat{A}.$$

We show later (in Corollary A.1, which follows this proof) that for every $s$ the set $\{P \in \Delta(\Omega) : s \in \hat{A}_P\}$ is Borel. Since

$$B \subset \{P \in \Delta(\Omega) : s \notin \hat{A}_P\}$$

for every $s \in \Omega - \hat{A}$, we obtain that

$$\forall_{s \in \Omega - \hat{A}} \qquad \zeta(\{P \in \Delta(\Omega) : s \in \hat{A}_P\}) = 0.$$

This means that $\Omega - \hat{A} \subset R_\zeta^0$, and so the complement of $R_\zeta^0$ is a first category set.

We will now construct sets $\hat{A}$ and $B$ with properties (A.1)–(A.2). Consider the sets

$$B_k^i(m) := \{P \in \Delta(\Omega) : t(i, k, P) > m\};$$

Lemma A.1, which follows this proof, shows that the sets $B_k^i(m)$ are open and so Borel. Since this sequence of sets is descending (with respect to $m$, for any given $k$ and $i$) and its intersection is empty, for every $l = 1, 2, \ldots$ there exists an $m$ such that

$$\zeta(B_k^i(m)) \leq \frac{1}{2^{k+i+l}};$$

denote by $m_k^i(l)$ any such $m$.

Let now

$$A_k^i(l) := \Omega - C(s_m^i) \qquad \text{for } m = m_k^i(l),$$

$$A_k(l) := \bigcap_{i=1}^{\infty} A_k^i(l)$$

and

$$\hat{A}(l) := \bigcup_{k=1}^{\infty} A_k(l) \cup \bigcup_{i=1}^{\infty} \{s^i\}.$$

The set $A(l)$ is a countable union of closed sets with empty interior by an argument analogous to that used, in the main body of the paper, for the case of the sets $\hat{A}_P$, and so is

$$\hat{A} := \bigcup_{l=1}^{\infty} \hat{A}(l).$$

To show that (A.1) and (A.2) are satisfied, notice that, by the definition of $B_k^i(m)$ for $m = m_k^i(l)$, if $P \notin B_k^i(m)$, then $C(s_m^i) \subset C(s_{t(i,k,P)}^i)$; therefore,

$$\text{if } P \in \Delta(\Omega) - \bigcup_{i=1}^{\infty} B_k^i(m_k^i(l)) \qquad \text{then } A_k(P) \subset A_k(l),$$

which in turn yields that

$$\text{if } P \in \Delta(\Omega) - \bigcup_{k=1}^{\infty} \bigcup_{i=1}^{\infty} B_k^i(m_k^i(l)) \qquad \text{then } \hat{A}_P \subset \hat{A}(l).$$

Thus,

$$B := \bigcup_{l=1}^{\infty} \left[ \Delta(\Omega) - \bigcup_{k=1}^{\infty} \bigcup_{i=1}^{\infty} B_k^i(m_k^i(l)) \right] \subset \{P \in \Delta(\Omega) : \hat{A}_P \subset \hat{A}\}.$$

It remains to show that $\zeta(B) = 1$; however,

$$\zeta(B) \geq \zeta\left(\Delta(\Omega) - \bigcup_{k=1}^{\infty} \bigcup_{i=1}^{\infty} B_k^i(m_k^i(l))\right)$$

$$\geq 1 - \sum_{k=1}^{\infty} \sum_{i=1}^{\infty} \zeta(B_k^i(m_k^i(l)))$$

$$\geq 1 - \sum_{k=1}^{\infty} \sum_{i=1}^{\infty} \frac{1}{2^{k+i+l}} = 1 - \frac{1}{2^l}$$

for every $l \in N$.  $\square$

LEMMA A.1.  *For every $t \in N$, the set*

$$B_k^i(t) := \{P \in \Delta(\Omega) : t(i,k,P) > t\}$$

*is open.*

PROOF.  Let $P \in B_k^i(t)$. By definition, $t(i,k,P) > t$, which means that

$$P(C(s_t^i) - \{s^i\}) > \frac{1}{2^{k+i}}.$$

Further, there exists an $m > t$ such that

$$(A.3) \qquad\qquad P(C(s_t^i) - C(s_m^i)) > \frac{1}{2^{k+i}};$$

indeed, the sequence of sets $C(s_t^i) - C(s_m^i)$ is ascending (as $m$ goes to infinity), and its union is equal to $C(s_t^i) - \{s^i\}$.

Note that each cylinder is an open and closed subset of $\Omega$, and so is $C(s_t^i) - C(s_m^i)$. Thus, the function $f : \Omega \to R$ given by

$$f(s) = \begin{cases} 1, & \text{for every } s \in C(s_t^i) - C(s_m^i), \\ 0, & \text{for every } s \notin C(s_t^i) - C(s_m^i), \end{cases}$$

is continuous.

Let

$$(A.4) \qquad\qquad \delta := P(C(s_t^i) - C(s_m^i)) - \frac{1}{2^{k+i}},$$

and let $N(P)$ stand for the set all measures $Q \in \Delta(\Omega)$ such that

$$\left| \int f \, dQ - \int f \, dP \right| < \delta.$$

This last inequality means that

$$|Q(C(s_t^i) - C(s_m^i)) - P(C(s_t^i) - C(s_m^i))| < \delta;$$

by (A.3) and (A.4),

$$Q(C(s_t^i) - \{s^i\}) \geq Q(C(s_t^i) - C(s_m^i)) > 1/2^{k+i},$$

which implies that $Q \in B_i^k(t)$.

That is, the set $N(P)$, which is an open neighborhood of $P$ in weak$^*$-topology, is contained in $B_k^i(t)$.   □

COROLLARY A.1.   *For every $s \in \Omega$, the set $\{P \in \Delta(\Omega) : s \in \hat{A}_P\}$ is Borel.*

PROOF.   If $s = s^i$ for some $i \in N$, then $\{P \in \Delta(\Omega) : s \in \hat{A}_P\} = \Delta(\Omega)$ by definition. Suppose, therefore, that $s \neq s^i$ for any $i \in N$. It suffices to show that the sets $\{P \in \Delta(\Omega) : s \in C(s_{t(i,k,P)}^i)\}$ are Borel, as

$$\{P \in \Delta(\Omega) : s \in \hat{A}_P\} = \Delta(\Omega) - \left[ \bigcap_{k=1}^{\infty} \left( \bigcup_{i=1}^{\infty} \{P \in \Delta(\Omega) : s \in C(s_{t(i,k,P)}^i)\} \right) \right].$$

Since $s \neq s^i$ and $C_P^k(s_0^i) = \Delta(\Omega)$, there is a unique $m = 0, 1, \ldots$ such that $s \in C_P^k(s_m^i) - C_P^k(s_{m+1}^i)$. For this $m$, we have that $\{P \in \Delta(\Omega) : s \in C(s_{t(i,k,P)}^i)\} = \Delta(\Omega) - \{P \in \Delta(\Omega) : t(i,k,P) > m\}$, and so the set $\{P \in \Delta(\Omega) : s \in C(s_{t(i,k,P)}^i)\}$ is closed by Lemma A.1.   □

PROOF OF PROPOSITION 2.   First, consider $P \in (\mathcal{A})^c$. It follows from (4.2) and from $P(R_P^m) \leq \frac{1}{(m+1)^3}$ that, for any cylinder $C \subseteq R_P^m$,

$$Q_P^m(C) \geq \frac{P(C)}{P(R_P^m)} \geq (m+1)^3 P(C).$$

By definition, $Q_P(A) \geq \pi(m) Q_P^m(A) = \frac{1}{(m+1)m} Q_P^m(A)$ for any $m \geq 1$ and set $A \in \Im$. Hence, for any cylinder $C \subseteq R_P^m$,

(A.5)      $$Q_P(C) \geq \frac{(m+1)^3}{(m+1)m} P(C) \geq m P(C).$$

Given a theory $P \in (\mathcal{A})^c$, let $\hat{R}_P$ and $\bar{R}_P$ be the rejection sets of $P$ for the tests $\hat{T}_S$ and $\bar{T}$, respectively. Consider a path $s \in \bar{R}_P$. Then, $s \in \hat{R}_P^m$ for every $m \geq 1$, because $\hat{T}^m$ is harder than $\hat{T}_S$. Since $\hat{R}_P^m$ comprises cylinders, there is a cylinder $C(s_{n(m)})$ with base on $s_{n(m)} = s \mid n(m)$ such that $C(s_{n(m)}) \subseteq \hat{R}_P^m \subseteq R_P^m$. By (A.5),

$$\frac{Q_P(C(s_{n(m)}))}{P(C(s_{n(m)}))} \geq m \quad \Longrightarrow \quad \lim_{m \to \infty} \frac{Q_P(C(s_{n(m)}))}{P(C(s_{n(m)}))} = \infty.$$

Thus, $s \in \bar{R}_P$. So, if $P \in (\mathcal{A})^c$, then $\hat{R}_P \subseteq \bar{R}_P$.

Now, consider a theory $P \in \mathcal{A}$. Consider a path $s \in \hat{R}_P$. Given that $P$ assigns probability one to the union of all its atoms, it follows that $s$ is not an atom of $P$. So, $P(C(s_t))$, $s_t = s \mid t$, approaches zero as $t$ goes to infinity. Since $P$ has finitely many atoms, there exists an $\eta > 0$ such that the probability of each atom is smaller than $\eta$. Let $t$ be large enough so that $P(C(s_t)) < \eta$. It follows that $C(s_t)$ contains no atom. Hence, $P(C(s_t)) = 0$. So, $s \in \bar{R}_P$. It follows that $\bar{T}$ is harder than $\hat{T}$. □

REMARK A.1. The existence of a cylinder with finite base that satisfies (4.1) can be shown as follows.

If $P \in (\mathcal{A})^c$, then there are infinitely many atoms, or the set of all nonatoms has positive measure. If there are infinitely many atoms, then one of them (say $s$) has to have measure below $\frac{1}{4(m+1)^3}$. It follows that the $P(C(s_t))$, $s_t = s \mid t$, must be strictly positive and (for $t$ large enough) smaller than $\frac{1}{2(m+1)^3}$. Now, consider the case in which the set of all nonatoms of $P$ has positive measure. Assume that every nonatom $s \in \Omega$ of $P$ is contained in a cylinder $C_s$ (with finite base) such that $P$ assigns zero probability to $C_s$. Let $\bar{C}$ be the union of all zero-probability cylinders. There are only countably many cylinders with finite base. So, $P(\bar{C}) = 0$. Given that $\bar{C}$ contains all nonatoms of $P$, it follows that the set of nonatoms has zero measure. This is a contradiction. So, $P$ has a nonatom $s$ such that $P(C(s_t)) \neq 0$, $s_t = s \mid t$, for all $t \in N$. Since $P(C(s_t))$ approaches zero as $t$ goes to infinity, $P(C(s_t))$ must (for $t$ large enough) be smaller than $\frac{1}{2(m+1)^3}$.

LEMMA A.2. *The test $\widetilde{T}^{LR}$ satisfies the joint measurability condition.*

PROOF. We need show that the set

$$\{(\theta, s, P) \in \Theta \times \Omega \times \Delta(\Omega) : T_S^{LR}(\theta, s, P) = 1\}$$

is measurable. Notice that

$$\{(\theta, s, P) \in \Theta \times \Omega \times \Delta(\Omega) : T_S^{LR}(\theta, s, P) = 1\}$$

$$= \Big\{(\theta, s, P) \in \Theta \times \Omega \times \Delta(\Omega) :$$

$$\forall_{n \in N} P(C(s_n)) \neq 0 \text{ and } \limsup_{n \to \infty} \frac{Q_{S_\theta}(C(s_n))}{P(C(s_n))} < \infty\Big\}$$

$$= \Theta \times \Omega \times \Delta(\Omega)$$

$$- \bigcap_{k=1}^{\infty} \bigcup_{n=1}^{\infty} \{(\theta, s, P) \in \Theta \times \Omega \times \Delta(\Omega) : Q_{S_\theta}(C(s_n)) > k P(C(s_n))\},$$

where (as always) $s_n = s \mid n$. To see the last equality, notice that $Q_{S_\theta}(C(s_n)) > 0$ for every cylinder $C(s_n)$, so if $P(C(s_n)) = 0$, then the inequality $Q_{S_\theta}(C(s_n)) > kP(C(s_n))$ is satisfied. Therefore, it suffices to show that sets of the form

$$\{(\theta, s, P) \in \Theta \times \Omega \times \Delta(\Omega) : Q_{S_\theta}(C(s_n)) > kP(C(s_n))\}$$

are measurable.

Further, the inequality that defines this last set depends only on $s_n$, not on the entire path $s$, so this set can be represented as a union of sets

$$\{(\theta, P) \in \Theta \times \Delta(\Omega) : Q_{S_\theta}(C(s_n)) > kP(C(s_n))\} \times C(s_n),$$

and so it suffices to show that sets of the form $\{(\theta, P) \in \Theta \times \Delta(\Omega) : Q_{S_\theta}(C) > kP(C)\}$, where $C$ is a given cylinder (with base on $s_n$), are measurable. We will show that every set of this form is open.

Indeed, a pair $(\theta, P)$ belonging to this set means that

$$(A.6) \qquad \sum_{i:s_\theta^i \in C} \frac{1}{i(i+1)} > kP(C);$$

denote by $\eta > 0$ the difference between the two expressions. Take the set of all paths $\overline{\theta}$ such that

$$(A.7) \qquad \sum_{i:s_{\overline\theta}^i \in C} \frac{1}{i(i+1)} > M := kP(C) + \frac{\eta}{2},$$

and the set of all probability distributions $\overline{P}$ such that

$$|P(C) - \overline{P}(C)| < \frac{\eta}{2k}.$$

Obviously, the Cartesian product of the two sets contains the pair $(\theta, P)$, and any pair $(\overline{\theta}, \overline{P})$ that belongs to this Cartesian product satisfies condition (A.6) (for $P$ replaced with $\overline{P}$ and $\theta$ replaced with $\overline{\theta}$). So, it suffices to show that the two sets are open. The fact that the latter set is open follows directly from the definition of weak*-topology, as any cylinder $C$ is a closed and open subset of $\Omega$.

We will now show that the former set is also open. If $\overline{\theta}$ satisfies (A.7), then

$$\sum_{i:s_{\overline\theta}^i \in C^t(\overline\theta)} \frac{1}{i(i+1)} > M,$$

for some $t = 1, 2, \ldots$, where

$$C^t(\overline{\theta}) := C \cap \{0,1\}^t \times \{I_{t+1}(\overline{\theta})\} \times \{I_{t+2}(\overline{\theta})\} \times \cdots.$$

This follows from the assumption that the paths from $\{0,1\}^t \times \{I_t(\theta)\} \times \{I_{t+1}(\theta)\} \times \cdots$ precede the paths from $S_\theta - \{0,1\}^t \times \{I_{t+1}(\theta)\} \times \{I_{t+2}(\theta)\} \times$ $\cdots$ in the sequence $s_\theta^1, s_\theta^2, \ldots$ and the fact that, if an infinite sum exceeds some number, then sufficiently large finite sums exceed that number as well.

With no loss of generality, one can assume that $t \geq n$; recall that $n$ is the length of the base of $C$. Then, for any $\overline{\theta}' \in C(\theta_t)$, the cylinder with base on the first $t$ outcomes of $\overline{\theta}$, we have

$$\sum_{i:\, s_{\overline{\theta}}^i \in C} \frac{1}{i(i+1)} > \sum_{i:\, s_{\overline{\theta}'}^i \in C^t(\overline{\theta}')} \frac{1}{i(i+1)} = \sum_{i:\, s_{\overline{\theta}}^i \in C^t(\overline{\theta})} \frac{1}{i(i+1)} > M;$$

the middle equality follows from the fact that the first $t$ outcomes of $\overline{\theta}$ and $\overline{\theta}'$ coincide, and so $s^i(\overline{\theta}) \in C^t(\overline{\theta})$ is equivalent to $s^i(\overline{\theta}') \in C^t(\overline{\theta}')$. $\quad \square$

LEMMA A.3. *Let $S = \{s^1, s^2, \ldots\}$ be a countable dense subset of $\Omega$. Take any $P \in \Delta(\Omega)$ such that $P(S) = 0$. The set of paths on which $T_S^{LR}$ rejects $P$ is a superset of the set of paths on which the global category test $\hat{T}_S$ rejects $P$.*

PROOF. Take any $s \in \Omega$ such that $\hat{T}_S(s, P) = 0$ and $k = 1, 2, \ldots$. Then, $s \in C(s_{t(i,k,P)}^i)$ for some $i = 1, 2, \ldots$; in other words, $s_{t(i,k,P)} = s_{t(i,k,P)}^i$. On one hand, as $P(S) = 0$,

$$P(C(s_{t(i,k,P)}^i)) = P(C(s_{t(i,k,P)}^i) - \{s^i\}) \leq \frac{1}{2^{k+i}},$$

and, on the other,

$$Q^S(C(s_{t(i,k,P)}^i)) \geq Q^S(s^i) = \frac{1}{i(i+1)}.$$

Thus,

$$\frac{Q^S(C(s_{t(i,k,P)}))}{P(C(s_{t(i,k,P)}))} \geq \frac{2^{k+i}}{i(i+1)} \geq 2^{k-1}$$

for every $k = 1, 2, \ldots$, which means that the sequence $(\frac{Q^S(C(s_t))}{P(C(s_t))})_{t=1}^\infty$ is unbounded. $\quad \square$

PROOF OF PROPOSITION 3. We show that, except countably many paths $\theta$, the set $R_\zeta^0(\theta)$ is a superset of the revelation set $R_\zeta^0$ of the global category test $\hat{T}_S$ for $S = S_\theta$. Suppose that $R_\zeta^0(\theta)$ is not a superset of the revelation set $R_\zeta^0$ of the global category test $\hat{T}_S$ for $S = S_\theta$. Then, by Lemma A.3, it must be the case that

$$\zeta(\{P \in \Delta(\Omega) : P(S_\theta) > 0\}) > 0.$$

If this is the case for an uncountable number of $\theta \in \Theta$, then there exist $m, n \in N$ such that

$$\zeta\left(\left\{P \in \Delta(\Omega) : P(S_\theta) > \frac{1}{m}\right\}\right) > \frac{1}{n},$$

for an uncountable set $\Xi \subset \Theta$.

Notice that the set $\Theta$ can be partitioned into countable subsets such that $\theta^1$ and $\theta^2$ belong to the same subset if $S_{\theta^1} = S_{\theta^2}$; that is, paths $\theta^1$ and $\theta^2$ coincide on all but a finite number of outcomes. Notice further that, if $\theta^1$ and $\theta^2$ belong to the distinct subsets, then the sets $S_{\theta^1}$ and $S_{\theta^2}$ are disjoint. Without loss of generality, we can assume that distinct paths from $\Xi$ belong to distinct subsets; that is, they differ on an infinite number of outcomes. Therefore, for any $P \in \Delta(\Omega)$, $P(S_\theta) > 1/m$ for at most $m - 1$ paths $\theta \in \Xi$. Further, we can restrict attention to an infinite but countable subset of $\Xi$; from now on, we will denote this subset by $\Xi$.

Consider any linear ordering $\precsim$ of the set $\Xi$, and define sets

$$D_\theta^k := \left\{P \in \Delta(\Omega) : \theta \text{ is } k\text{th path from } \Xi \text{ such that } P(S_\theta) > \frac{1}{m}\right\}$$

for $k = 1, \ldots, m - 1$ and $\theta \in \Xi$.

For a given $k = 1, \ldots, m - 1$, the sets $D_\theta^k$ are pairwise disjoint, and, for a given $\theta \in \Xi$, the sets $D_\theta^k$ are pairwise disjoint. The measurability of the sets $D_\theta^k$ follows from the measurability of the set of all measures $P$ that at a given path $s$ have an atom of measure larger than a given number. Finally,

$$\left\{P \in \Delta(\Omega) : P(S_\theta) > \frac{1}{m}\right\} = D_\theta^1 \cup \cdots \cup D_\theta^{m-1}$$

for every $\theta \in \Xi$.

Thus,

$$m - 1 = (m - 1) \cdot \zeta(\Delta(\Omega))$$

$$\geq \sum_{k=1}^{m-1} \zeta\left(\bigcup_{\theta \in \Xi} D_\theta^k\right) = \sum_{k=1}^{m-1} \left(\sum_{\theta \in \Xi} \zeta(D_\theta^k)\right)$$

$$= \sum_{\theta \in \Xi} \left(\sum_{k=1}^{m-1} \zeta(D_\theta^k)\right)$$

$$= \sum_{\theta \in \Xi} \zeta\left(\left\{P \in \Delta(\Omega) : P(S_\theta) > \frac{1}{m}\right\}\right)$$

$$> \sum_{\theta \in \Xi} \frac{1}{n} = \infty,$$

a contradiction. $\square$

PROOF OF PROPOSITION 4.  Let $\tilde{T}' : \Omega \times \Delta(\Omega) \to \{0, 1\}$ be a test defined by $\tilde{T}'(s, P) = 1$ if and only if $\tilde{T}(\theta, s, P) = 1$ $\tilde{v}$-almost surely. The test $\tilde{T}'$ is a prequential test because, given two theories $P$ and $P'$ equivalent along $s$, for every $\theta \in \Theta$, $\tilde{T}(\theta, s, P) = 1$ if and only if $\tilde{T}(\theta, s, P') = 1$. Hence, $\tilde{T}(\theta, s, P) = 1$, $\tilde{v}$-almost surely, if and only if $\tilde{T}(\theta, s, P') = 1$, $\tilde{v}$-almost surely. Since, for every $\theta \in \Theta$,

$$P(\{s \in \Omega \mid \tilde{T}(\theta, s, P) = 1\}) = 1,$$

by Fubini's theorem

$$P(\{s \in \Omega \mid \tilde{T}(\theta, s, P) = 1 \ \tilde{v}\text{-almost surely}\}) = 1.$$

So, $\tilde{T}'$ does not reject the data-generating process with probability one. By Shmaya's (2008) result [which relies on Martin's (1998) theorem], there exists $\tilde{\zeta}$ such that for all $s \in \Omega$

$$\tilde{\zeta}(\{P \in \Delta(\Omega) \mid \tilde{T}'(s, P) = 1\}) \geq 1 - \delta. \qquad \square$$

Let $X$ be a metric space. Recall that a function $f : X \to R$ is *lower semicontinuous* at an $x \in X$ if, for every sequence $(x_n)_{n=1}^{\infty}$ converging to $x$,

$$\forall_{\varepsilon > 0} \ \exists_N \ \forall_{n \geq N} \qquad f(x_n) > f(x) - \varepsilon.$$

The function $f$ is lower semi-continuous if it is lower semi-continuous at every $x \in X$. We refer the reader to Engelking [(1989), Problem 1.7.14] for these definitions and some basic results regarding lower semi-continuous functions.

LEMMA A.4.  *Let $U \subset X$ be an open set where $X$ is a compact metric space. Equip $X$ with the $\sigma$-algebra of Borel subsets. Let $\Delta(X)$ be the set of all probability measures on $X$. Equip $\Delta(X)$ with the weak\*-topology. The function $F : \Delta(X) \to [0, 1]$ defined by*

$$F(P) = P(U)$$

*is lower semi-continuous.*

PROOF.  See Dudley (1989), Theorem 11.1.1(b).  $\square$

THEOREM [Fan (1953)].  *Let $X$ be a compact Hausdorff space, which is a convex subset of a linear space, and let $Y$ be a convex subset of linear space (not necessarily topologized). Let $f$ be a real-valued function on $X \times Y$ such that for every $y \in Y$, $f(x, y)$ is lower semi-continuous, with respect to $x$. If $f$ is also convex, with respect to $x$, and concave, with respect to $y$, then*

$$\min_{x \in X} \sup_{y \in Y} f(x, y) = \sup_{y \in Y} \min_{x \in X} f(x, y).$$

We note that Fan's (1953) theorem allows for $X$ and $Y$ that may not be subsets of linear spaces. We, however, apply his result only to subsets of linear spaces.

PROOF OF PROPOSITION 5. Let $X = \Delta(\Omega)$, let $Y$ be the subset of $\Delta(\Delta(\Omega))$ that consists of all random generators of theories with finite support. So, an element $\zeta$ of $Y$ can be described by a finite sequence of probability measures $\{P_1, \ldots, P_n\}$ and positive weights $\{\pi_1, \ldots, \pi_n\}$ that add up to one, where $\zeta$ selects $P_i$ with probability $\pi_i$, $i = 1, \ldots, n$. Let the function $f : X \times Y \to R$ be defined by

$$(A.8) \qquad f(P, \zeta) := E^P E^\zeta T = \sum_{i=1}^{n} \pi_i \int T(s, P_i) \, dP(s).$$

We now check that the assumptions of Fan's theorem are satisfied. Since $T$ is an open test, the set

$$U_Q = \{s \in \Omega : T(s, Q) = 1\}$$

is open for every $Q \in \Delta(\Omega)$. Therefore, by Lemma A.4,

$$P(U_Q) = \int T(s, Q) \, dP(s)$$

is a lower semi-continuous function of $P$. Thus, for every $\zeta \in Y$, the function $f(P, \zeta)$ is lower semi-continuous on $X$ as a weighted average of lower semi-continuous functions.

By definition, $f$ is linear with respect to both $x$ and $y$, and so it is convex with respect to $x$ and concave with respect to $y$. By the Riesz and Banach–Alaoglu theorems, $X$ is a compact space in weak*-topology; it is a metric space, and so Hausdorff [see, e.g., Rudin (1973), Theorem 3.17].

Thus, by Theorem 1,

$$\min_{P \in X} \sup_{\zeta \in Y} E^P E^\zeta T = \sup_{\zeta \in Y} \min_{P \in X} E^P E^\zeta T.$$

Notice that the left-hand side of this equality exceeds $1 - \varepsilon$, as the test $T$ is assumed not to reject the data-generating process with probability $1 - \varepsilon$; indeed, for a given $P \in X$, take $\zeta$ such that $\zeta(\{P\}) = 1$. Therefore, the right-hand side exceeds $1 - \varepsilon$, which yields the existence of a random generator of theories $\zeta \in Y$ such that

$$E^P E^\zeta T > 1 - \varepsilon - \delta$$

for every $P \in \Delta(\Omega)$. Taking, for any given $s \in \Omega$, the probability measure $P$ such that $P(\{s\}) = 1$, we obtain

$$\zeta(\{Q \in \Delta(\Omega) : T(s, Q) = 1\}) > 1 - \varepsilon - \delta. \qquad \qquad \square$$

**Acknowledgments.** We are extremely grateful to Peter Grünwald for pointing out an error in a previous version of this paper, suggesting the tests in Sections 4.1 and 4.2, stating and proving several results regarding their properties and encouraging the discussion in Sections 3.1 and 3.2. We thank Nabil Al-Najjar, Eddie Dekel, Yossi Feinberg, Dean Foster and Rakesh Vohra for their comments. We also thank workshop participants at Brown, Caltech, Carnegie Melon, Columbia, Cornell, Harvard, MIT, IMPA, Northwestern, Penn State, Princeton, Stanford, Rice, Tucson, U. of Chicago, U. of Bonn, U. of British Columbia, Berkeley, UCLA, University of Texas at Austin, University of Toronto, Washington University at St. Louis, Yale and at the Summer Meetings of the Econometric Society in Minneapolis 2006.

Department of Economics
Northwestern University
2001 Sheridan Road
Evanston, Illinois 60208
USA
E-mail: wo@northwestern.edu

Department of Economics
University of Pennsylvania
3718 Locust Walk
Philadelphia, Pennsylvania 19104
USA
E-mail: sandroni@sas.upenn.edu
  and
Kellogg School of Management
MEDS Department
2001 Sheridan Road
Evanston, Illinois 60208
USA
E-mail: sandroni@kellogg.northwestern.edu